\documentclass[titlepage,11pt]{article}
\usepackage{latexsym}
\usepackage{amsfonts}
\usepackage{amsmath}
\usepackage{mathtools}
\usepackage{xcolor}
\usepackage{comment}
\newcommand\blackslug{\hbox{\hskip 1pt \vrule width 4pt height 8pt depth 1.5pt
        \hskip 1pt}}
\newcommand\bbox{\hfill \quad \blackslug \bigbreak}
\newcommand\ind{\operatorname{ind}}

\def\LL{,\ldots,}
\newcommand{\vare}{\varepsilon}
\newcommand{\cupcup}{\cup \cdots\cup}
\def\poly{\operatorname{poly}}
\DeclarePairedDelimiter\abs{\lvert}{\rvert}%
\DeclarePairedDelimiter\ceil{\lceil}{\rceil}%
\DeclarePairedDelimiter\floor{\lfloor}{\rfloor}%
\title{Induced subgraph density. I. A loglog step towards Erd\H{o}s-Hajnal}
\author{
Matija Buci\'c\thanks{
Supported by the Oswald Veblen Fund and NSF grant CCF-1900460.}\\
Princeton University and \\Institute for Advanced Study,\\ Princeton, NJ 08544, USA
\and 
Tung Nguyen\thanks{Supported by AFOSR grant  FA9550-22-1-0234 and by NSF grant  DMS-2154169.}\\
Princeton University,\\ Princeton, NJ 08544, USA
\and 
Alex Scott\thanks{Research supported by EPSRC grant EP/X013642/1.}\\
Mathematical Institute,\\ University of Oxford, 
\\
Oxford OX2 6GG, UK
\and
Paul Seymour\footnotemark[2]\\
Princeton University,\\ Princeton, NJ 08544, USA}

\date{September 17, 2022; revised \today}

\newtheorem{thm}{}[section]

\newcommand{\Proof}{\noindent{\bf Proof.}\ \ }

\begin{document}
\maketitle
\begin{abstract}
In 1977, Erd\H{o}s and Hajnal made the conjecture that, for every graph $H$, there exists $c>0$ such that 
every $H$-free graph $G$ has a clique or stable set of size at least $|G|^c$;
and they proved that this is true with $ |G|^c$ replaced by $2^{c\sqrt{\log |G|}}$.
Until now, there has been no improvement on this result (for general $H$).

We prove a strengthening: that for every graph $H$, there exists $c>0$ such that
every $H$-free graph $G$ with $|G|\ge 2$ has a clique or stable set of size at least 
$$2^{c\sqrt{\log |G|\log\log|G|}}.$$
Indeed, we prove the corresponding strengthening of a theorem of Fox and Sudakov, which in turn was a common strengthening of theorems
of R\"odl, Nikiforov, and the theorem of Erd\H{o}s and Hajnal mentioned above.

\end{abstract}

\section{Introduction}

A graph $G$ {\em contains} a graph $H$ if $H$ is isomorphic to an induced subgraph of $G$, and $G$ is {\em $H$-free} otherwise.
$|G|$ denotes the number of vertices of the graph $G$; and we write $\kappa(G)$ for the largest~$t$ such that $G$ has a 
clique or stable set of cardinality $t$.
For most $n$-vertex graphs $G$, $\kappa(G)= O(\log n)$, but this changes dramatically if we forbid some induced subgraph.
In 1977, Erd\H{o}s and Hajnal~\cite{EH0,EH} proposed the following well-known conjecture:

\begin{thm}\label{EHconj}
{\bf Conjecture: }For every graph $H$ there exists $c>0$ such that $\kappa(G)\ge |G|^c$ for every $H$-free graph $G$.
\end{thm}

The Erd\H{o}s-Hajnal conjecture has attracted a great deal of attention over the years, but despite this, it is only known to be true for a few 
graphs $H$. Until recently, it was only known for the graphs with at most five vertices except the five-vertex path and its complement, and graphs that can be 
made from these by vertex-substitution. In two recent papers~\cite{density4, density7}, three of us have got further: we have shown it for the five-vertex 
path, and for infinitely many other 
graphs that are cannot be built from smaller graphs by vertex-substitution. But it remains the case that 
the graphs known to satisfy \ref{EHconj} are rare and highly restricted.

On the other hand, it is known that excluding any fixed induced subgraph will guarantee that $\kappa(G)$ is much bigger than $\log |G|$.
Erd\H{o}s and Hajnal~\cite{EH} themselves proved the following:

\begin{thm}\label{EHthm}
For every graph $H$ there exists $c>0$ such that $\kappa(G)\ge 2^{c\sqrt{\log |G|}}$ for every non-null $H$-free graph $G$.
\end{thm}
Indeed they proved something slightly stronger, that for all $H$ and all $c>0$ the same conclusion holds, provided that $|G|$ is 
sufficiently large.
Rather surprisingly, until now there has been no improvement on this for a general graph $H$.
Our result is such an improvement:
\begin{thm}\label{mainthm}
For every graph $H$ there exists $c>0$ such that $\kappa(G)\ge 2^{c\sqrt{\log |G|\log\log |G|}}$ for every $H$-free graph $G$ 
with $|G|\ge 2$.
\end{thm}

Over the years, there have been several theorems discovered that are related to the Erd\H{o}s-Hajnal conjecture, 
and our proof method allows us to strengthen some of them.
First, a fundamental theorem of R\"odl~\cite{rodl} shows that every $H$-free graph contains a large dense or sparse set:
\begin{thm}\label{rodl}
For every graph $H$ and all $x>0$, there exists $\delta>0$ with the following property.
For every $H$-free graph $G$, there exists $S\subseteq V(G)$ with $|S|\ge \delta|G|$ such that 
one of $G[S],\overline G[S]$ has  at most $x\binom{|S|}{2}$ edges.
\end{thm}

How large can $\delta$ be as a function of $x$?
Fox and Sudakov~\cite{foxsudakov} conjectured that $\delta$ can be taken to be a polynomial in $x$, and noted that this would imply the Erd\H os-Hajnal conjecture.
R\"odl's original proof gave a tower-type bound, because it used the regularity
lemma, but Fox and Sudakov~\cite{foxsudakov} made a significant improvement, proving a version of \ref{rodl} with better bounds:
\begin{thm}\label{foxsudakov}
There exists $c>0$ such that for every graph $H$ and all $x\in (0,1/2)$, \ref{rodl} holds with
$\delta= 2^{-c|H|(\log \frac{1}{x})^2}$.
\end{thm}
This result implies \ref{EHthm} (by setting $x=2^{-\sqrt{\frac{\log |G|}{2c|H|}}}$, and applying Tur\'an's theorem;
the proof is similar to the proof that
\ref{mainthmnik} implies \ref{mainthm}, which we give later).

Nikiforov~\cite{nikiforov} gave a different strengthening of \ref{rodl}, allowing for a small number of copies of $H$:
\begin{thm}\label{nikiforov}
For every graph $H$ and all $x>0$, there exists $\delta>0$ such that if $G$ is a graph 
containing fewer than $(\delta|G|)^{|H|}$ induced copies of $H$,
then there exists $S\subseteq V(G)$ with $|S|\ge \delta|G|$
such that one of $G[S], \overline{G}[S]$ has at most $x\binom{|S|}{2}$ edges.
\end{thm}

Fox and Sudakov~\cite{foxsudakov} were able to incorporate the analogous strengthening of \ref{nikiforov} into \ref{foxsudakov}:
\begin{thm}\label{foxsudakovnik}
There exists $c>0$ such that for every graph $H$ and all $x\in (0,1/2)$, setting
$\delta= 2^{-c|H|(\log \frac{1}{x})^2}$
satisfies \ref{nikiforov}.
\end{thm}

Our main result is: 
\begin{thm}\label{mainthmnik}
For every graph $H$ there exists $c$ such that, if $x \in (0,1/2)$ 
and 
\[\delta=2^{-c(\log\frac1x)^2/\log\log\frac1x},\]
%$$\delta= \exp{\left(\frac{-c\left(\log \frac{1}{\vare}\right)^2}{\log\log \frac{1}{\vare}}\right)},$$
and $G$ is a graph containing fewer than $(\delta|G|)^{|H|}$ induced copies of $H$,
then there exists $S\subseteq V(G)$ with $|S|\ge \delta|G|$
such that one of $G[S], \overline{G}[S]$ has at most $x\binom{|S|}{2}$ edges.
\end{thm}
\ref{mainthmnik} strengthens the result \ref{foxsudakovnik} of Fox and Sudakov, and improves the best known quantitative bounds in
Nikiforov's theorem \ref{nikiforov} and R\"odl's theorem \ref{rodl}.
It also implies \ref{mainthm},  as we will show later. 
The proof of \ref{mainthmnik} is by induction on $|H|$: let us mention that it is essential for inductive purposes that 
we allow $G$ to contain a few copies of $H$, rather than that $G$ is $H$-free.

The rest of the paper is organized as follows.
We give some definitions and outline the strategy in section \ref{sec:b}.  We prove \ref{tunglemma}
in section \ref{sec:T} and \ref{manycopies} in section \ref{sec:M}.  We complete the proofs of our main results in section \ref{sec:D}, where we 
show that \ref{manycopies} implies \ref{mainthmnik}, and that \ref{mainthmnik} implies \ref{mainthm}.
In the final section we discuss analogous results for tournaments and ordered graphs.

We note that all logarithms in this paper are to base 2.

%%%%%%%%%%%%%%%%%%%%%%%%%%%%%%%%%%%%%%%%%%%%%%%%%%%%%%%%%%%%%%%%%%%%%%%%%%%%%%%%%%%

\section{Blockades, and a sketch of the proof}\label{sec:b}

A {\em cograph} means a $P_4$-free graph, where $P_4$ denotes the path with four vertices; and we denote by $\mu(G)$
the largest $t$ such that some $t$-vertex induced subgraph of $G$ is a cograph. Cliques and stable sets induce cographs, and
every cograph $J$ has a clique or stable set of
size at least $|J|^{1/2}$; so \ref{EHconj}, \ref{EHthm}, and \ref{mainthm} are equivalent to the same statements with
$\kappa(G)$ replaced by $\mu(G)$, and in that form they are often easier to work with.

If $X\subseteq V(G)$, $G[X]$ denotes the induced subgraph with vertex set $X$.
A {\em pure pair} in $G$ is a pair of disjoint subsets $A,B$ of $V(G)$ such that either there are no edges between $A,B$, or
all edges between $A,B$ are
present.
If we are trying to prove that
$\mu(G)\ge f(|G|)$ for all $H$-free graphs $G$, where $f$ is some function, it is enough to know that all $H$-free graphs $G$
with $|G|>1$ have pure pairs $A,B$ with $|A|,|B|$ appropriately large in terms of $|G|$. Because then we could
deduce by induction
on $|G|$ that $G[A]$ contains a cograph $C$ with $|C|\ge f(|A|)$, and similarly $G[B]$ contains a large cograph $D$, and so
$V(C)\cup V(D)$ induces a cograph in $G$, and therefore $\mu(G)\ge f(|A|)+f(|B|)$; and if $|A|,|B|$ are large enough, then
$f(|A|)+f(|B|)\ge f(|G|)$ and the inductive step is complete.  The key factor here is how large we can take $A$ and $B$ to be.  For example, if we can take them to be a constant fraction of $G$ then we would obtain a bound of form $|G|^c$ (see, for example, \cite{pure1, pure2}).

In fact, for this purpose, we do not really need the pair $A,B$ to be pure.
Suppose that either every vertex in $B$ has at most  $|A|/(2\mu(G))$ neighbours in $A$,
or every vertex in $B$ has at most  $|A|/(2\mu(G))$ non-neighbours in $A$.
Choose $D\subseteq B$ as before; then, since $|D|\le \mu(G)$, there
exists $A'\subseteq A$ with $|A'|\ge |A|/2$ such that $A',D$ is a pure pair, and we apply the inductive hypothesis to $G[A']$,
and
reach the same conclusion as before.

Do such pairs $A,B$ necessarily exist with $A,B$ large? If $A,B\subseteq V(G)$ are disjoint, we say $B$ is {\em $x$-sparse to $A$} in $G$ if
every vertex in $B$ has at most  $x|A|$ neighbours in $A$. 
Erd\H{o}s and Hajnal~\cite{EH} proved:
\begin{thm}\label{EHstart}
For every graph $H$, there exists $c>0$
such that for every $H$-free graph $G$ with $|G|\ge 2$, and all $x\in (0,1/2)$, there exist
disjoint $A,B\subseteq V(G)$ with $|A|,|B|\ge x^{c} |G|$, such that 
$B$ is $x$-sparse to $A$ in one of $G,\overline{G}$.
\end{thm}
Then they used \ref{EHstart} (with an appropriate choice of $x$)  and the inductive argument sketched above, to prove \ref{EHthm}.
But perhaps \ref{EHstart} can be strengthened. There is a pretty conjecture of Conlon, Fox and Sudakov~\cite{conlon} that would strengthen it:
\begin{thm}\label{conlon}
{\bf Conjecture:} For every graph $H$ there exists $c_1,c_2> 0$ such that for every $H$-free graph $G$ with $|G|\ge 2$, and all
$x\in (0,1/2)$, there exist
disjoint $A,B\subseteq V(G)$ with $|A|\ge x^{c_1} |G|$ and $|B|\ge c_2|G|$, such that 
$B$ is $x$-sparse to $A$ in one of $G,\overline{G}$.
\end{thm}
If this were true, then, as Conlon, Fox and Sudakov observed, the inductive argument would yield exactly our result \ref{mainthm} (by choosing $x =1/(2\mu(G))$).

Our argument takes a different approach and leaves this conjecture open, however.
Rather than look for a sparse pair of large subsets, we will look for a large number of smaller subsets, with a sparseness condition between each pair. 

If $G,H$ are graphs, a {\em copy} of $H$ in $G$ is an isomorphism from $H$ to an induced
subgraph of $G$, and we denote by $\ind_H(G)$ the number of copies of $H$ in $G$. 
A {\em blockade} in $G$ is a sequence $\mathcal{B}=(B_1\LL B_k)$ of
pairwise disjoint
subsets of $V(G)$, and we call $B_1\LL B_k$ its {\em blocks}. (In some earlier papers, the blocks of a blockade must be nonempty,
but here it is convenient to allow empty blocks.)
The {\em length} of the blockade $\mathcal{B}=(B_1\LL B_k)$ is $k$, and its {\em width} is the minimum of the cardinalities of its blocks. For $\vare>0$, the blockade
$\mathcal{B}=(B_1\LL B_k)$ is
{\em $x$-sparse} in $G$ if for all $i$ with $1\le i\le k$, $B_{i+1}\cupcup B_k$ is $\vare$-sparse to $B_i$ in $G$; and 
{\em $\vare$-restricted} if for all $i$ with $1\le i\le k$, $B_{i+1}\cupcup B_k$ is $\vare$-sparse to $B_i$ in one of $G,\overline{G}$.

Let us give a sketch of the proof of \ref{mainthm}. The key step is the following lemma, which says that if $G$ does not contain too many copies of $H$ then we can find a large blockade that is very dense or very sparse:
\begin{thm}\label{manycopies}
For all $H$, there exist $k_1,k_2>0$ such that for every non-null graph $G$ and every $x$ with
$0<x\le\frac{1}{8\vert H\vert}$,
if $\ind_H(G)<
x^{k_1}|G|^{|H|}$, then there is an $x$-restricted blockade in $G$ of length at least $2\log(1/x)$, and
width at least 
$\lfloor x^{k_2}|G|\rfloor$.
\end{thm}

First let us see that \ref{manycopies} implies \ref{mainthm}. Choose $c>0$ sufficiently small; we will show that $c$ satisfies \ref{mainthm} by 
induction on $|G|$.  Let $x:=1/(2\mu(G))$. It is easy to arrange that $x\le\frac{1}{8\vert H\vert}$, and so we can apply \ref{manycopies} to the $H$-free
graph $G$,
and obtain a blockade $(B_1\LL B_k)$. Then choose subsets $D_i\subseteq B_i$ for $i = k,k-1\LL 1$ in turn, such that $D_i\cupcup D_{k}$ induces a cograph, as follows.
Having chosen $D_{i+1}\LL D_k$, since $D_{i+1}\cupcup D_k$ induces a cograph, it has cardinality at most $\mu(G)$, and so
(assuming every vertex of $B_{i+1}\cupcup B_k$ has at most $x|B_i|$
neighbours in $B_i$; the other case is similar), at least half the vertices in $B_i$ have no neighbour in $D_{i+1}\cupcup D_k$.
By induction,
we may choose $D_i$ from this half, inducing a cograph and
of cardinality at least $(|B_i|/2)^c$. Then $D_{i}\cupcup D_k$ induces a cograph, completing the inductive definition.
Consequently $\mu(G)\ge \sum (|B_i|/2)^c$,
and the result follows after some calculation,
which we omit. (We will not actually prove \ref{mainthm} this way; our proof goes via the stronger theorem \ref{mainthmnik}, which can also be derived
from \ref{manycopies} with more work.)

The main issue is how to prove \ref{manycopies}.
The core of the argument is the following lemma, which applies to graphs with few copies of $H$ and has two possible outcomes: either we can drop to a large induced subgraph with very few copies of some subgraph $H'$ with $|H'|<|H|$, in which case we can continue by induction; or we obtain two large sets of vertices that have very few edges or very few nonedges between them, in which case we attempt to build a blockade.

\begin{thm}\label{tunglemma}
Let $H$ be a graph, let $g\in V(H)$, and let $H':= H\setminus \{g\}$. Let $b,c> 0$, and let $a:=b+(1+c)|H|$.
Let $G$ be a graph, let
$A,B$ be disjoint subsets of $V(G)$, and let $0<x\le 1/2$. Suppose that every vertex in $A$ has at least
$x|B|$ non-neighbours in $B$. Then either:
\begin{itemize}
\item there exists $B'\subseteq B$ with $|B'|\ge x|B|$ such that $\ind_{H'}(G[B'])< x^b|B'|^{|H'|}$; or
\item $\ind_H(G)\ge x^{a}|A|\cdot|B|^{|H|-1}$; or
\item there exists $A'\subseteq A$ and $B'\subseteq B$ with $|A'|\ge x^{a}|A|$ and $|B'|\ge x^{a}|B|$
such that the number of edges between $A',B'$ is at most $2x^c|A'|\cdot |B'|$.
\end{itemize}
\end{thm}
In the remainder of this section, we will sketch a proof of \ref{tunglemma}, and then sketch how we use it to prove \ref{manycopies}.

\medskip\noindent {\em Sketch of proof of \ref{tunglemma}:}
The idea of the proof of \ref{tunglemma} is as follows. Let $g$ have degree $d$ in $H$, let $H_d:=H$, and for $i=d-1,d-2\LL 0$
let $H_i$ be obtained from $H_{i+1}$ by deleting one of the $i+1$ edges of $H_{i+1}$ incident with $g$. We are interested in copies of $H_i$ in $G$,
where $g$ is mapped into $A$ and the other vertices of $H_i$ are mapped into $B$. (Let us call such copies ``special''.) Each $v\in A$
has at least $x|B|$ non-neighbours in $B$, and we may assume that this set of non-neighbours ($B'$ say) induces a subgraph that contains at least
$x^{b}|B'|^{|H|-1}\ge x^{b+|H|}|B|^{|H|-1}$ copies of $H\setminus \{g\}$, because otherwise the first bullet holds. Since this is true for each $v$,
there are
at least $x^{b+|H|}|A|\cdot |B|^{|H|-1}$ special copies of $H_0$. On the other hand, we may assume that there are fewer than $x^{a}|A|\cdot|B|^{|H|-1}$ 
special copies of $H_d=H$, because otherwise the second bullet holds. 
So for some $t$, the number of special copies of $H_t$ is at least $x^{b+|H|+ct}|A|\cdot |B|^{|H|-1}$, and the number of special copies of $H_{t+1}$ is less than
$x^{b+|H|+c(t+1)}|A|\cdot |B|^{|H|-1}$. Let us focus on this value of $t$.

Now $H_t$ was obtained from $H_{t+1}$ by deleting some edge incident with $g$, say $gh$.  
Fix a copy $\phi$ of $H\setminus \{g,h\}$ in $G[B]$; let $U$ be the set of $u\in B$ such that mapping
$h$ to $u$ extends $\phi$ to a copy of $H\setminus \{g\}$, and let $V$ be the set of $v\in A$ such that mapping $g$ to $v$ extends $\phi$ to a copy of 
$H\setminus \{h\}$. The number of special copies of $H_t$ is the sum, over all $\phi$, of the number of nonedges between $U,V$, and the number of special
copies of $H_{t+1}$ is the sum over all $\phi$ of the number of edges between $U,V$. Since the second sum is at most $x^c$ times the first, 
there is a choice
of $\phi$ such that the number of edges between $U,V$ is at most $x^c|U|\cdot |V|$ (and by allowing a factor of two here and averaging, we can arrange that also  
$|U|\cdot |V|\ge x^a|A|\cdot|B|$). But then $U,V$ satisfy the third bullet, and this will prove \ref{tunglemma}.

\medskip\noindent {\em Sketch of proof of \ref{manycopies}:}
Next, let us explain how to use \ref{tunglemma} to prove \ref{manycopies}.
Let $g\in V(H)$ and $F=H\setminus \{g\}$. We will assume inductively that
\ref{manycopies} holds for $F$, with $k_1,k_2$ replaced by $k_1', k_2'$ say.
Choose $k_1,k_2$ sufficiently large, and 
let $G$ be a graph with $\ind_H(G)<
x^{k_1}|G|^{|H|}$. We must show that there is an $x$-restricted blockade in $G$ of length at least $\log(1/x)$, and
width at least
$\lfloor x^{k_2}|G|\rfloor$.

Choose an induced subgraph $J$ of $H$ maximal such that
$G$ contains a large ``approximate blowup'' of $J$; that is, $|J|$ disjoint subsets $A_j\; (j\in V(J))$ of $V(G)$,
each of size about $x^k|G|$ (where $k$ is an appropriate constant depending on $|J|$),
and such that for all distinct $i,j\in V(J)$, if $ij\notin E(J)$ 
then $A_i,A_j$ are $x$-sparse to each other in $G$, and the same in the complement if $ij\in E(J)$.
It cannot be that $J=H$ since otherwise there would be $x^{k_1}|G|^{|H|}$ copies of $H$ in $G$,
contrary to the hypothesis; let $h\in V(H)\setminus V(J)$, and let $J':=H[V(J)\cup \{h\}]$. There is no large approximate blowup of $J'$ in $G$ (even allowing
its 
sets to be a little smaller than the $A_j$'s, and the sparsity between them to be relaxed a little), and we will exploit this. Let $W$ be the set of vertices of $G$ in none of the sets $A_j\;(j\in V(J))$.
Thus $W$ contains almost all vertices of $G$.

Let us assume, first, that there exists $X_0\subseteq W$ with $|X_0|\ge |W|/2$, such that for each $j\in V(J)$, if $hj\in E(H)$ then every vertex in $X_0$ has at least $x|A_j|$
neighbours in $A_j$, and if $hj\notin E(H)$ then every vertex in $X_0$ has at least $x|A_j|$ non-neighbours in $A_j$. Then we can obtain a blockade with 
the properties we want, as follows. Let $j\in V(J)$, and suppose that $h,j$ are 
nonadjacent in $H$ (the other case is the same in the complement). Each vertex in $X_0$ has at least $x|A_j|$ non-neighbours in $A_j$, and so we can apply
\ref{tunglemma} with $X_0,A_j$ in place of $A,B$. If the first outcome of \ref{tunglemma} holds, then since $F=H\setminus \{g\}$ satisfies \ref{manycopies}, 
there is an $x$-restricted blockade in $G[A_j]$ of length at least $\log(1/x)$, and
width at least
$\lfloor x^{k_2'}|A_j|\rfloor$; and that blockade has the properties we want, if we arrange the constants properly. The second outcome
of \ref{tunglemma} cannot hold, since it would contradict the hypothesis of \ref{manycopies}. So we assume that the third outcome holds; and we can therefore choose
$X_1\subseteq X_0$ and $C_j\subseteq A_j$ such that $|X_1|\ge \poly(x)|X_0|$, and $|C_j|\ge \poly(x) |A_j|$, such that there are at 
most $\poly(x)|X_1|\cdot|A_j'|$ edges between $X_1,C_j$. 
Note what has happened: we started with every vertex in $X_0$ just having a few (at least $x|A_j|$) non-neighbours in $A_j$, and now there are almost no edges 
between $X_1$ and $C_{j}$. We can assume that $|C_j|$ has size about equal to its lower bound $\poly(x) |A_j|$;
and by removing some ``outlier'' vertices from $X_1$, we can assume that in addition, 
$X_1$ is $\poly(x)$-sparse to $C_j$. (The advantage of this is that, we will choose successively smaller subsets of $X_1$, and they will all be $\poly(x)$-sparse to $C_j$.)

Now choose some other vertex $j'\in V(J)$ different from $j$, and apply \ref{tunglemma} to the pair $X_1,A_{j'}$; we may assume this gives us
$X_2\subseteq X_1$ with $|X_2|\ge \poly(x)|X_1|$, and $C_{j'}\subseteq A_{j'}$ with $|C_{j'}|\ge \poly(x)|A_{j'}|$, such that 
$X_2$ is $\poly(x)$-sparse to $C_{j'}$ in $G$ if $ij'\notin E(H)$, and the same in the complement if $ij'\in E(H)$.
Continue in this way until we have processed each vertex of $J$. This gives us $X_{|J|}\subseteq W$ with $|X_{|J|}|\ge \poly(x)|W|$, 
and $C_j\subseteq A_j$ with $|C_j|\ge \poly(x) |A_j|$ for each $j\in V(J)$, such that
$X_{|J|}$ is $\poly(x)$-sparse to $C_j$ in $\overline{G}$ or in $G$ (depending whether $hj\in E(H)$ or not) for each $j\in V(J)$. 
We can assume that each $C_j$ is $\poly(x)$-sparse to $X_{|J|}$ (again, by removing
a few outliers). 
Since the sets $A_j\;(j\in V(J))$ are $\poly(x)$-sparse to each other in $G$ or $\overline{G}$ (where $\poly(x)$ is some polynomial of large degree), as in the definition of an approximate
blowup, and since $|C_j|\ge \poly(x)|A_j|$ for each $j$, it follows that the sets $C_j\;(j\in V(J))$ are still $\poly(x)$-sparse to each other in $G$ or 
$\overline{G}$ (where $\poly(x)$ is now some polynomial of somewhat smaller degree).
This gives an 
approximate blowup of $J':=H[V(J)\cup \{h\}]$, 
contradicting the choice of $J$. So this cannot happen; and therefore, at some step, the first outcome of \ref{tunglemma} held, and so we obtained the blockade we want.

Consequently, we may assume that there is no such $X_0$; and so, for some $j\in V(J)$, there is a subset $A\subseteq W$ with $|A|\ge |W|/(2|H|)$ that is  $x$-sparse to $A_j$
in $G$ or in $\overline{G}$. 
Now repeat the proof, working
completely within $A$. After $2\log(1/x)$ iterations of the argument, we will produce
the $x$-restricted blockade that we want; and until that stage, the various subsets we must deal with are still large enough that the argument is valid.
This completes the sketch of the proof of \ref{mainthm}.

%%%%%%%%%%%%%%%%%%%%%%%%%%%%%%%%%%%%%%%%%%%%%%%%%%%%%%%%%%%%%%%%%%%%%%%%%%%%%%%%%%%%%%%%%%%%%%%%%%%%%%%%%%%%%%%%%%%%%%%%%%%%%

\section{The proof of \ref{tunglemma}}\label{sec:T}

In this section we prove \ref{tunglemma}, which we restate:

\begin{thm}\label{tunglemma2}
Let $H$ be a graph, and let $g\in V(H)$. Let $b,c> 0$, and define $a:=b+(1+c)|H|$.
Let $G$ be a graph, let
$A,B$ be disjoint subsets of $V(G)$, and let $0<x\le 1/2$, such that $A$ is $(1-x)$-sparse to $B$.
Then either:
\begin{itemize}
\item there exists $B'\subseteq B$ with $|B'|\ge x|B|$ such that $\ind_{H\setminus \{g\}}(G[B'])< x^b|B'|^{|H|-1}$; or
\item $\ind_H(G)\ge x^{a}|A|\cdot|B|^{|H|-1}$; or
\item there exists $A'\subseteq A$ and $B'\subseteq B$ with $|A'|\ge x^{a}|A|$ and $|B'|\ge x^{a}|B|$
such that the number of edges between $A',B'$ is at most $2x^c|A'|\cdot |B'|$.
\end{itemize}
\end{thm}
\Proof 
Let $g$ have degree $d$, let $H_d:=H$, and inductively for $t =d-1\LL 0$, let $H_t$ be obtained from $H_{t+1}$
by deleting one of the $t+1$ edges of $H_{t+1}$ incident with $g$.  Let $k:=|H|$, and $s:=|A|\cdot |B|^{k-1}$.

Let $v\in A$. By hypothesis, the set $B'$ (say) of non-neighbours of $v$ in $B$ has cardinality at least 
$x|B|$. If $G[B'] $
contains fewer than  $x^b|B'|^{k-1}$ copies of $H\setminus \{g\}=H_0\setminus \{g\}$, then the first bullet holds, so we assume not. Consequently there are at least 
$x^b|B'|^{k-1}\ge x^{k-1+b}|B|^{k-1}$ 
copies $\phi$ of $H_0$
in $G$ such that $\phi(g)=v$ and $\phi(j)\in B$ for each $j\in V(H)\setminus \{g\}$. It follows by 
summing over all $v\in A$ 
that there are at least
$x^{k-1+b}s$
copies $\phi$ of $H_0$ such that $\phi(g)\in A$ and $\phi(j)\in B$ for each $j\in V(H)\setminus \{g\}$.

For $0\le t\le d$, let $\tau_t$ be the number of copies $\phi$ of $H_t$ in $G$ such that $\phi(g)\in A$ and $\phi(j)\in B$ 
for each $j\in V(H)\setminus \{g\}$.
We have just seen that  $\tau_0\ge x^{k-1+b}s$.
We may assume that $\tau_d<sx^{a}\le x^{k-1+b+cd}s$,
because otherwise the second bullet holds.
Consequently, for some $t$ with $1\le t\le d$, 
\begin{itemize}
\item $\tau_{t-1}\ge x^{k-1+b+c(t-1)}s\ge 2x^as$, and
\item $\tau_{t}< x^{k-1+b+ct}s$, and therefore $\tau_t< x^c\tau_{t-1}$.
\end{itemize}
Let $\Phi$ be the set of all copies $\phi$ of $H_{t-1}$ in $G$ such that $\phi(g)\in A$ and $\phi(j)\in B$ for each $j\in V(H)\setminus \{g\}$.
There is one edge of $H_t$ that is not an edge of $H_{t-1}$, say $gh$. Let $J$ be the graph obtained from $H$ by deleting both $g$ 
and $h$, and let $\Psi$ be the set of all copies of $J$ in $G[B]$. Each member of $\Phi$ is an extension of some member of $\Psi$.
For each $\psi\in \Psi$, let $n(\psi)$ be the number of $\phi\in \Phi$ that are extensions of $\psi$.
Thus the sum of $n(\psi)$ over all $\psi\in \Psi$ equals
$\tau_{t-1}$. Let $\Psi'$ be the set of all $\psi\in \Psi$ such that $n(\psi)\ge \tau_{t-1}/(2|B|^{k-2})$, and 
$\Psi''=\Psi\setminus \Psi'$. Thus
$$\sum_{\phi\in \Psi''}n(\psi)\le |B|^{k-2}\left(\frac{\tau_{t-1}}{2|B|^{k-2}}\right)= \tau_{t-1}/2$$
since $| \Psi''|\le |B|^{k-2}$ and $n(\psi)\le \tau_{t-1}/(2|B|^{k-2}) $ for each $\phi\in \Psi''$.
Since 
$$\sum_{\psi\in \Psi''}n(\psi) + \sum_{\psi\in \Psi'}n(\psi)= \tau_{t-1},$$
it follows that 
$$\sum_{\psi\in \Psi'}n(\psi)\ge \tau_{t-1}/2.$$
For each $\psi\in \Psi'$, let $U(\psi)$ be the set of all $u\in B$ such that mapping $i$ to $u$ extends $\psi$ to a copy of $H\setminus \{g\}$,
and let $V(\psi)$ be the set of all $v\in A$ such that mapping $g$ to $v$ extends $\psi$ to a copy of $H\setminus \{h\}$. Thus
$n(\psi)$ is the number of pairs $(u,v)$ with $u\in U(\psi)$ and $v\in V(\psi)$ such that $u,v$ are nonadjacent. Let $p(\psi)$
be the number of edges between $U(\psi)$ and $V(\psi)$. Thus $\tau_t$ is at least the sum of $p(\psi)$ over all $\psi\in \Psi'$.
Since $\tau_t<x^c\tau_{t-1}$, it follows that
$$\sum_{\psi\in \Psi'}p(\psi)\le \tau_t<x^c\tau_{t-1}\le 2x^c\sum_{\psi\in \Psi'}n(\psi);$$
and consequently there exists $\psi\in \Psi'$ such that $p(\psi)\le 2x^cn(\psi)$.

Since $\psi\in \Psi'$, it follows that 
$$|U(\psi)|\cdot |V(\psi)|\ge n(\psi)\ge \frac{\tau_{t-1}}{2|B|^{k-2}}\ge \frac{2sx^{a}}{2|B|^{k-2}}
=x^{a}|A|\cdot |B|;$$ 
and since $|U(\psi)|\le |A|$ and $|V(\psi)|\le |B|$, it follows that 
$|U(\psi)|\ge x^{a}|A|$
and similarly $|V(\psi)| \ge x^{a}|B|$. Since 
$$p(\psi)\le 2x^cn(\psi)\le 2x^c|U(\psi)|\cdot |V(\psi)|,$$
there are at most $2x^c|U(\psi)|\cdot |V(\psi)|$ edges between $U(\psi)$ and $V(\psi)$. Hence the third bullet holds, setting $A'=U(\psi)$
and $B'=V(\psi)$. This proves \ref{tunglemma2}.~\bbox

%%%%%%%%%%%%%%%%%%%%%%%%%%%%%%%%%%%%%%%%%%%%%%%%%%%%%%%%%%%%%%%%%%%%%%%%%%%%%

\section{The proof of \ref{manycopies}}\label{sec:M}
Now we turn to the proof of \ref{manycopies}. We will need:

\begin{thm}\label{tidy}
Let $A,B$ be disjoint subsets of $V(G)$, such that there are at most $c|A|\cdot |B|$ edges between $A$ and $B$. Then there 
exists $A'\subseteq A$ with $|A'|\ge |A|/2$ such that $A'$ is $2c$-sparse to $B$.
\end{thm}
\Proof
There are at most $c|A|\cdot |B|$ edges between $A$ and $B$, and so at most $|A|/2$ vertices in $A$ have more than
$2c|B|$ neighbours in $B$. This proves \ref{tidy}.~\bbox

Let $J$ be a graph, and $t>0$ an integer, and $0\le q\le 1$ a real number. Let $G$ be a graph, and let $A_j\;(j\in V(J))$
be pairwise disjoint subsets of $V(G)$. We say that the family $(A_j:j\in V(J))$
is a {\em $(t,q)$-blowup of $J$} if
\begin{itemize}
\item each set $A_j\;(j\in V(J))$ has cardinality $t$;
\item for all distinct $i,j\in V(J)$, if $ij\notin E(J)$ then $A_i, A_j$ are $q$-sparse to each other in $G$,
and if 
$ij\in E(J)$ then $A_i, A_j$ are $q$-sparse to each other in $\overline{G}$.
\end{itemize}

We observe:
\begin{thm}\label{blowuptocopies}
Let $J$ be a graph, and $t>0$ an integer. If there is a $(t,1/|J|)$-blowup of $J$ in $G$, then 
$\ind_J(G)\ge (t/|J|)^{|J|}$.
\end{thm}
\Proof
Let $(A_j:j\in V(J))$ be a $(t,1/|J|)$-blowup of $J$ in $G$. If $I$ is an induced subgraph of $J$, a copy $\phi$ of $I$ is {\em good}
if $\phi(i)\in A_i$ for each $i\in I$. 
\\
\\
(1) {\em Let $I$ be an induced subgraph of $J$, 
and suppose that $\phi$ is a good copy of $I$. Then there are at least $(t/|J|)^{|J|-|I|}$ good copies of $J$ that extend $\phi$.}
\\
\\
The proof is by induction on $|V(J)|-|V(I)|$. If this is zero then the claim is true, so we may assume that 
there exists $j\in V(J)\setminus V(I)$. Let $I'$ be the induced subgraph of $J$ with vertex set $V(I)\cup \{j\}$.
For each $i\in V(I)$, let us say that $v\in A_j$ is {\em $i$-conforming}
if either $ij\in E(J)$ and $\phi(i),v$ are adjacent in $G$, or $ij\notin E(J)$ and $\phi(i),v$ are nonadjacent in $G$.
From the definition of 
a $(t,q)$-blowup, for each $i\in V(I)$ there are at most $t/|J|$ vertices in $A_j$ that are not $i$-conforming; and so there are 
at least $t-t|I|/|J|\ge t/|J|$ vertices $v\in A_j$ such that $v$ is $i$-conforming for each $i\in V(I)$. For each such $v$, 
let $\phi'$ be the extension of $\phi$ obtained by mapping $j$ to $v$; then $\phi'$ is a good copy of $I'$. From the inductive 
hypothesis, there are at least $(t/|J|)^{|J|-|I|-1}$ good copies of $J$ that extend $\phi'$; and since there are at least
$t/|J|$ choices of $v$ and hence of $\phi'$, the claim follows. This proves (1).

\bigskip

But then the theorem follows from (1) by setting $I$ to be the null graph. This proves \ref{blowuptocopies}.~\bbox

The bulk of the proof of \ref{manycopies} consists of the following lemma:

\begin{thm}\label{manycopies0}
For all graphs $H$, all $g\in V(H)$, and all $\alpha>0$, there exist $\beta,\gamma >0$ such that for every graph $G$ with $|G|\ge 2$
and all
$x$ with $0<x\le 1/(8|H|)$,
either:
\begin{itemize}
\item there exists $A\subseteq V(G)$ with $|A|\ge x^{\beta}|G|$ such that $\ind_{H\setminus \{g\}}(G[A])<
x^{\alpha}|A|^{|H|-1}$; or
\item $\ind_H(G)\ge x^{\gamma}|G|^{|H|}$; or
\item there are disjoint subsets $A,B\subseteq V(G)$ with $|A|\ge x^{\beta}|G|$ and $|B|\ge |G|/(2|H|)$, such that $B$ is 
$x$-sparse to $A$ in one of $G,\overline{G}$.
\end{itemize}
\end{thm}

\Proof We may assume that $|H|\ge 2$, because otherwise the theorem holds taking $\gamma=1$.
It suffices to prove the result assuming that $1/x$ is an integer. Indeed, suppose that $H,g,\alpha$ are given, and
setting $\beta=\beta'$ and $\gamma =\gamma'$ satisfies the theorem for all $G$ and $x$ with $1/x$ an integer.
Then setting $\beta=2\beta'$ and $\gamma =2\gamma'$ satisfies the theorem for all $G$ and $x$. To see this, let
$0<x\le 1/(8|H|)$, and let $x'=1/\left(\lceil 1/x \rceil\right)$. Then $1/x'$ is an integer, and 
$1/x'= \lceil 1/x \rceil\le \frac{8|H|+1}{8|H|x}$, and so 
$$x^2\le\frac{x}{8|H|}\le  \frac{8|H|x}{8|H|+1} \le x'\le x.$$ 
Consequently
$(x')^{\beta'} \ge  (x^2)^{\beta'} = x^\beta$ and similarly $(x')^{\gamma'} \ge x^\gamma$
and hence, whichever bullet of the theorem holds for $x', \beta',\gamma'$, the same bullet holds for $x,2\beta',2\gamma'$. So
to prove the theorem, 
we just need to exhibit values of $\beta,\gamma$ that work when $1/x$ is an integer.

By increasing $\alpha$ if necessary, we may assume that $\alpha$ is an integer and $\alpha\ge |H|(|H|+1)$.
Define $r_{|H|}=0$, and inductively for $i = |H|-1\LL 1$ define
$$r_i:= \alpha+2|H|+1 +(|H|+1)r_{i+1}.$$
Let 
$\beta:=r_1+3$
%\beta-r_k \ge 3
%\beta\ge 2+k r_k-(k-1)r_{k+1}
% 0 \le \beta-2-r{k+1}-k r_k+kr_{k+1}
and $\gamma:= 2r_1+ \beta|H|$.
%$\gamma\ge (\beta+1)|H|$
%\gamma\ge 2r_1+ \beta|H|
We claim that
$\beta, \gamma$ satisfy the theorem (when $1/x$ is an integer).

Thus, let $G$ be a graph with $|G|\ge 2$,  and let $x>0$ such that $0<x\le 1/(8|H|)$, where $1/x$ is an integer.
If $x^{\beta}|G|\le  1$, the third bullet is true taking $|A|=1$ (unless $|G|-1<|G|/|H|$, which is impossible since $|G|,|H|\ge 2$),
so we
may assume that $x^{\beta}|G|> 1$.
We assume the first two bullets of the theorem are false, and we will show that the third holds.

Let $t:=\lfloor x^{\beta-1}|G|\rfloor $; thus $t\ge x^{\beta}|G|$, because
$x^{\beta-1} |G|\ge 1$
and $x\le 1/2$. Let $t_i:= x^{-r_i}t$ and $q_i:=x^{r_i}/|H|$ for $1\le i\le |H|$.
Thus $t_1\LL t_{|H|}$ are integers.

Since $\gamma/|H|\ge \beta+1$, it follows that 
$$t\ge x^{\beta}|G|\ge (1/x) x^{\gamma/|H|}|G|\ge |H| x^{\gamma/|H|}|G|,$$ 
and consequently
%%%%%%%%%%%$\gamma/|H|\ge \beta+1)$
$(t/|H|)^{|H|}\ge x^{\gamma}|G|^{|H|}$.
Hence by \ref{blowuptocopies}, there is no $(t,1/|H|)$-blowup (that is, no $(t_{|H|},q_{|H|})$-blowup) of $H$ in $G$.
Let $J$ be a maximal induced
subgraph of $H$ such that there is a $(t_{|J|},q_{|J|})$-blowup $(A_j:j\in V(J))$ of $J$ in $G$, and let $k:=|J|$.

Thus $J\ne H$; let $h\in V(H)\setminus V(J)$, and $L:=\bigcup_{j\in V(J)}A_j$. For each $j\in V(J)$, let $M_j$ be the set of
vertices $v\in V(G)\setminus L$ such that
\begin{itemize}
\item if $hj\in E(H)$, then $v$ has at most $x|A_j|$ neighbours in $A_j$;
\item if $hj\notin E(H)$, then $v$ has at most $x|A_j|$ non-neighbours in $A_j$.
\end{itemize}
For each $j\in V(J)$, since $t_k\ge t\ge x^{\beta}|G|$,  we may assume that $|M_j|< |G|/(2|H|)$,
since otherwise the third bullet of the
theorem holds. Since
$$|L|=kt_k\le kx^{\beta-1-r_k}|G|\le x^2|H|\cdot |G|\le |G|/(2|H|),$$
%%%%%%%%%%%%\beta-r_k \ge 3
it follows that the union of $L$ and the sets $M_j\;(j\in V(J))$ has
cardinality at most $|G|/2$. Let $Z$ be the set of vertices of $G$ that do not belong to $L$ or to any of the sets $M_j\;(j\in V(J))$.
Thus $|Z|\ge |G|/2$; and for each $j\in V(J)$, if
$hj\notin E(H)$ then $Z$ is $(1-x)$-sparse to $A_j$, and if $hj\in E(H)$ then $Z$ is $(1-x)$-sparse to $A_j$ in $\overline{G}$. Let
$s:=x^{r_k-r_{k+1}}$. Thus $t_{k+1} = t_k$ and $q_{k+1}=q_k/s$.
\\
\\
(1) {\em Let $j\in V(J)$, and let $Y\subseteq Z$ with $|Y|\ge |Z|s^{k-1}$. Then
there exist $C\subseteq A_j$ with $|C|=2t_{k+1}$, and $X\subseteq Y$ with $|X|\ge s|Y|$,
such that $X$ is $\frac12 q_{k+1}$-sparse to $C$ in $G$ if $hj\notin E(H)$, and $X$ is $\frac12 q_{k+1}$-sparse to $C$ in $\overline{G}$ if
$hj\in E(H)$.
}
\\
\\
By taking complements if necessary, we may assume that $hj\notin E(H)$, and so $Y$ is $(1-x)$-sparse to $A_j$.
We will apply \ref{tunglemma} with $b,c,A,B$ replaced by $\alpha,r_{k+1}+1,Y,A_j$; note that the expression
$b+(1+c)|H|$ of \ref{tunglemma} becomes $\alpha+(r_{k+1}+2)|H| = r_k-r_{k+1}-1$.
By \ref{tunglemma}, we deduce that either:
\begin{itemize}
\item there exists $A'\subseteq A_j$ with $|A'|\ge x|A_j|$ such that 
$\ind_{H\setminus \{g\}}H(G[A'])
<x^{\alpha}|A'|^{|H|-1}$; or
\item $\ind_H(G)\ge x^{r_k-r_{k+1}-1}|Y|\cdot|A_j|^{|H|-1}$; or
\item there exist $A'\subseteq A_j$ and $D\subseteq Y$ with $|A'|\ge x^{r_k-r_{k+1}-1}|A_j|$ and $|D|\ge x^{r_k-r_{k+1}-1}|Y|$
such that the number of edges between $A',D$ is at most $2x^{r_{k+1}+1}|A'|\cdot |D|$.
\end{itemize}

If the first bullet above holds, then the first bullet of the theorem holds, since
$|A'|\ge x|A_j|=xt_{|J|}=x^{1-r_{|J|}}t\ge t\ge  x^{\beta}|G|$
(because $|J|<|H|$ and so $r_{|J|}\ge 1$), a contradiction.
If the second holds, then the second bullet of the theorem holds, also a contradiction, since
$$x^{r_k}|Y|\cdot|A_j|^{|H|-1} \ge x^{r_k+ r_1+ \beta|H|}|G|^{|H|}\ge  x^{2r_1+ \beta|H|}|G|^{|H|}= x^{\gamma}|G|^{|H|}$$
(because 
$$|A_j|^{|H|-1}\ge t^{|H|-1}\ge x^{\beta(|H|-1)}|G|^{|H|-1} \ge x^{\beta|H|}|G|^{|H|-1}.$$
and $|Y|\ge |Z|s^{k-1}\ge x^{(k-1)r_k}|G|\ge x^{r_1}|G|$.)
%%%%%%%%%%%%%%%%%%%%\gamma\ge 2r_1+ \beta|H|

Thus the third bullet above holds. Let $A',D$ be the corresponding subsets. Since
$$|A'| \ge x^{r_k-r_{k+1}-1}t_{k}=x^{r_k-r_{k+1}-1}x^{-r_k}t = t_{k+1}/x\ge 2t_{k+1},$$
it follows by averaging that
we may choose $C\subseteq A'$ with $|C|=2t_{k+1}$
such that the number of edges between $C,D$ is at most $2x^{r_{k+1}+1}|C|\cdot |D|$.
By \ref{tidy}, there exists $X\subseteq D$ with
$$|X|\ge \frac{|D|}{2}\ge \frac12  x^{r_k-r_{k+1}-1}|Y|\ge s|Y|$$
such that $X$ is $4x^{r_{k+1}+1}$-sparse to $C$, and hence $\frac12 q_{k+1}$-sparse to $C$, since
%4x^{r_{k+1}}\le \frac12 q_{k+1}$
%ie (1/|H|) x^{r_{k+1}} \ge 8x^{r_{k+1}+1}$
%%%%%%%%%%%%%%%%%ie x^{|H|+1} \le 1/(8|H|)$
$4x^{r_{k+1}+1}\le \frac{1}{2|H|} x^{r_{k+1}}$ (because $x\le 1/(8|H|)$). This proves (1).

\bigskip

Starting with $Y=Z$, and applying (1) recursively to each $j\in V(J)$, we obtain a subset $X\subseteq Z$ with
$|X|\ge |Z|s^{k}$, and a subset $C_j\subseteq A_j$ with $|C_j|=2t_{k+1}$ for each $j\in V(J)$,
such that for each $j\in V(J)$, $X$ is $\frac12 q_{k+1}$-sparse to $C_j$ in $G$ if $hj\notin E(H)$, and $X$ is 
$\frac12 q_{k+1}$-sparse 
to $C_j$ in $\overline{G}$ if $hj\in E(H)$. Since (from the choice of $\beta$)
$$|X|\ge s^{k}|Z|\ge  \frac12 x^{k r_k-kr_{k+1}}|G|\ge x^{\beta-1-r_{k+1}}|G| \ge tx^{-r_{k+1}}=t_{k+1}$$
%%%%%%%%%%%%%%%% 0 \le \beta-2-r{k+1}-k r_k+kr_{k+1}
we may choose $D_h\subseteq X$ with $|D_h|=t_{k+1}$.
By \ref{tidy}, for each $j\in V(J)$ there exists $D_j\subseteq C_j$ with $|D_j|=t_{k+1}$
such that $D_j, D_h$ are $q_{k+1}$-sparse to each other in $G$ if $hj\notin E(H)$, and $D_j, D_h$ are $q_{k+1}$-sparse
to each other in $\overline{G}$ if $hj\in E(H)$. Let $i,j\in V(J)$ be distinct. We assume that $ij\notin E(H)$ (the other case is the same in the complement).
Since $(A_j:j\in V(J))$ is a $(t_k,q_k)$-blowup of $J$, it follows
that $A_i$ is $q_k$-sparse to $A_j$. 
Since $|D_j|=t_{k+1}=x^{-r_{k+1}}t=s|A_j|$, this implies that $A_i$ (and hence $D_i$) is $(q_k/s)$-sparse to $D_j$, that is,
$q_{k+1}$-sparse to $D_j$. Consequently 
$(D_j:j\in V(J'))$ is a $(t_{k+1},q_{k+1})$-blowup of $J'$, where $J'$ is the induced
subgraph of $H$ with vertex set $V(J)\cup \{h\}$,
contrary to the maximality of $J$. This proves \ref{manycopies0}.~\bbox

Now we use \ref{manycopies0} to prove \ref{manycopies}, which we restate, slightly strengthened:

\begin{thm}\label{manycopies2}
For all $H$, there exist $k_1,k_2>0$ such that for every non-null graph $G$ and every $x$ with 
%$0<x\le 1/(8|H|)$ and 
%$x\le e^{-c_1\sqrt{ \log |G|}}$, 
$0< x\le\frac{1}{8|H|}$,
if $\ind_H(G)< x^{k_1}|G|^{|H|}$, 
there is an $x$-restricted blockade in $G$ with length at least $2\log(1/x)$ and width at least  $\floor{x^{k_2}|G|}$.
Consequently, for all such $G,x$, there is a blockade in $G$ with length at least $\log(1/x)$ and width at least  $\floor{x^{k_2}|G|}$ that is $x$-sparse 
in one of $G,\overline{G}$.
\end{thm}
%$\beta+d\le k_2$
%\gamma+d|H|\le k_1
%$(c_2')^2 (1-(\beta+d)/k_2) \ge c_2^2$
%k_2\ge \beta+d$.
\Proof
The first statement is trivially true when $|H|\le 1$, and 
we proceed by induction on $|H|$. So 
we may assume that $|H|\ge 2$, and $g\in V(H)$, and the theorem holds for $H\setminus \{g\}$; let
$k_1',k_2'$ be the corresponding constants (using $H\setminus \{g\}$ instead of $H$).
Let $d:= \log (2|H|)$.
Choose $\beta,\gamma$ satisfying \ref{manycopies0}, 
taking $\alpha=k_1'$.
Let $k_1:=\gamma+2d|H|$ and $k_2:= k_2'+\beta+2d$. We will show that $k_1,k_2$ satisfy the theorem.

Thus, let $G$ be a graph and $x$ with $0<x\le\frac{1}{8|H|}$ 
%and $x\le e^{-c_1\sqrt{ \log |G|}}$, 
such that $\ind_H(G)< x^{k_1}|G|^{|H|}$.
Choose an $x$-restricted blockade $\mathcal{B}=(B_1\LL B_k)$ in $G$ with $k$ maximum such that $B_1\LL B_{k-1}$ have cardinality at least 
$x^{k_2}|G|$ and $|B_k|\ge (2|H|)^{1-k}|G|$. 
\\
\\
(1) {\em We may assume that $|B_k|> x^{2d}|G|$.}
\\
\\
We may assume that 
$k-1<2\log(1/x)$
and so
$$|B_k|\ge (2|H|)^{1-k}|G|=2^{(1-k)d}|G|> 2^{-2d\log (1/x)}|G| = x^{2d}|G|.$$
This proves (1).

\bigskip
If $x^{k_2}|G|< 1$, the result holds trivially (because $\floor{x^{k_2}\abs G}=0$ and blockades may contain empty blocks),
so we may assume that $|G|\ge x^{-k_2}$.
By (1), $\abs{B_k}\ge x^{2d}\abs G\ge x^{2d-k_2}>1$ and so $|B_k|\ge 2$.
Let us apply \ref{manycopies0} to $G[B_k]$, taking $\alpha=k_1'$. We deduce that
either:
\begin{itemize}
\item there exists $A\subseteq B_k$ with $|A|\ge x^{\beta}|B_k|$ such that 
$\ind_{H\setminus \{g\}}(G[A])<
x^{\alpha}|A|^{|H|-1}$; or
\item $\ind_H(G[B_k])\ge x^{\gamma}|B_k|^{|H|}$; or
\item there are disjoint subsets $A,B\subseteq B_k$ with $|A|\ge x^{\beta}|B_k|$ and $|B|\ge |B_k|/(2|H|)$, such that $B$ is 
$x$-sparse to $A$ in one of $G,\overline{G}$.
\end{itemize}
The second is impossible, since by (1), $x^{\gamma}|B_k|^{|H|}\ge x^{\gamma} x^{2d|H|} |G|^{|H|} = x^{k_1}|G|^{|H|}$.
Also, the third is impossible, from the maximality of $k$, because $x^\beta|B_k|\ge x^{\beta} x^{2d}|G|\ge x^{k_2}|G|$.
%$\beta+d\le k_2$
%\gamma+d|H|\le k_1
Thus the first holds. Let $A$ be the corresponding subset. 
Since
%$x\le e^{-c_1\sqrt{ \log |G|}}\le e^{-c_1\sqrt{ \log |B_k|}}$,
$\vert A\vert \ge x^{\beta}\vert B_k\vert \ge x^{\beta+d}\vert G\vert$,
the inductive hypothesis gives an $x$-restricted blockade in $G[A]$ with length at least 
$2\log(1/x)$ and width at least  $\floor{x^{k_2'}|A|}\ge \floor{x^{k_2'}x^{\beta+d}\vert G\vert}= \floor{x^{k_2}\vert G\vert}$.
This proves the first statement of the theorem.

For the second statement, let $G,x$ be as before, and let $(B_1\LL B_k)$ be an $x$-restricted blockade in $G$ with length at least $2\log(1/x)$ and width at least  
$\floor{x^{k_2}|G|}$. Let $I$ be the set of $i\in \{1\LL k\}$ such that $B_{i+1}\cupcup B_k$ is $x$-sparse to $B_i$ in $G$, and $J=\{1\LL k\}\setminus I$;
then for all $i\in J$, $B_{i+1}\cupcup B_k$ is $x$-sparse to $B_i$ in $\overline{G}$. So the blockade $(B_i:i\in I)$ is $x$-sparse in $G$, and 
$(B_i:i\in J)$ is $x$-sparse in $\overline{G}$, and one of them has length at least $k/2\ge \log(1/x)$. 
This proves \ref{manycopies2}.~\bbox

%%%%%%%%%%%%%%%%%%%%%%%%%%%%%%%%%%%%%%%%%%%%%%%%%%%%%%%%%%%%%%%%%%%%%%%%%%%%%%%%%%%%%%%%%%

\section{Deriving the main theorem}\label{sec:D}

It remains to show that \ref{manycopies} implies \ref{mainthmnik}, and that \ref{mainthmnik} implies \ref{mainthm},
and we do so in this section. 
\ref{manycopies} says that graphs that do not contain many copies of $H$ admit blockades with certain properties, but the length 
of this blockade is critical. \ref{manycopies} gives 
blockades of length $\log(1/x)$, but one might hope that for some graphs $H$, we could obtain a version of \ref{manycopies} 
that gave longer blockades; and then there would be corresponding improvements in \ref{mainthmnik} and \ref{mainthm}. 
With that in mind, we have written the argument of this section in greater generality than is needed for this paper.
(We will use this generality in~\cite{density4,density5}, for instance.) Let us say the {\em edge-density} of a nonnull graph $J$ is the number of edges of $J$
divided by $\binom{|J|}{2}$.

A function $\ell\colon(0,\frac12)\to\mathbb R^+$ is {\em subreciprocal} if it is non-increasing and $1<\ell(x)\le 1/x$ for all $x\in(0,\frac12)$.
(To prove \ref{mainthmnik} we will only need the subreciprocal function $\ell(x):=\log(1/x)$.)
If $\ell$ is a subreciprocal function, a graph $H$ is {\em $\ell$-divisive}
if there exist $c \in (0,\frac12)$ and $d > 1$ such that for every $x \in (0,c)$ and every graph $G$ with 
$\ind_H(G)\le x^{d}\abs G^{\abs H}$, there is a blockade $(B_1,\ldots,B_k)$ in $G$ with length at least $\ell(x)$ and width at least 
$\floor{x^{d}\abs G}$, that is $x$-sparse in one of $G,\overline{G}$.

Erd\H{o}s and Hajnal~\cite{EH} proved that every graph is $\ell$-divisive where $\ell(x)=2$ for $x\in (0,1/2)$; and \ref{manycopies} implies the following:
\begin{thm}
	\label{thm:extensive}
	Every graph is $\ell$-divisive where $\ell(x):=\log(1/x)$ for $0<x<1/2$.
\end{thm}

	The next theorem implies our main result \ref{mainthmnik}, by defining $\ell$ as in \ref{thm:extensive}.
	The proof is an adaptation of an argument of Fox and Sudakov~\cite{foxsudakov}.
\begin{thm}
	\label{thm:mainthm2}
Let $H$ be an $\ell$-divisive graph for some subreciprocal function $\ell$.
Then there exists $C>0$ such that, if $\varepsilon \in (0,\frac12 )$ and 
        \[\delta=2^{-C\log^2(1/\vare)/\log(\ell(\vare))},\]
        then for every graph $G$ with  $\ind_H(G)\le(\delta\abs G)^{\abs H}$, there exists
$S\subseteq V(G)$ with $\abs S\ge\delta\abs G$ such that one of $G[S],\overline G[S]$ has edge-density at most $\varepsilon$.
\end{thm}
\Proof
	Let $c\in(0,\frac12)$ and $d>1$, as in the definition of $\ell$-divisive. 
	Let $z:=\ell(c)^{-1/2}\in (0,1)$, 
%OLD $z:=\ell(c)^{-1/4}\in (0,1)$,
and let $b>2$ be such that $2^{2-b}=1-z$.
	We will first show that setting $C=20bd$ satisfies the theorem when $\varepsilon\in(0,c)$,
and then at the end of the proof, give a value of $C$ that works in general.

Thus, let  $\varepsilon\in(0,c)$, and choose $\delta$ such that 
$$\log(1/\delta)= \frac{20bd\log^2(1/\varepsilon)}{\log(\ell(\varepsilon))}.$$
	%	In what follows we assume $\varepsilon\in(0,\frac1he^{-64})$.
	%	Let $k_1,k_2>1$ be as in \cref{thm:blocks}. 
	Let $x:=\frac{1-z}{2}\varepsilon=2^{1-b}\vare$, and let $p:=\ell(x)z$;
%OLD $p:=z\sqrt{\ell(x)}$;
	then $p^2\ge \ell(x)>1$ 
%OLD $p\ge \sqrt[4]{\ell(x)}>1$
since $\ell(x)\ge\ell(\vare)\ge\ell(c)=z^{-2}$.
%OLD $\ell(x)\ge\ell(\vare)\ge\ell(c)=z^{-4}$.

	Let $t$ be the least integer such that $p^t\ge \varepsilon^{-2}$;
	then since $\ell(\vare)\le\min(\ell(x),1/\vare)
	\le\min(p^2,1/\vare)$, we obtain
%OLD \le\min(p^4,1/\vare)$,
	\[1\le t=\ceil*{\frac{2\log(1/\varepsilon)}{\log p}}
	\le\ceil*{\frac{4\log(1/\vare)}{\log(\ell(\vare))}}
	\le\frac{5\log(1/\varepsilon)}{\log(\ell(\varepsilon))}.\] 
	Let $\eta:=\frac14x^d$;
	then $x=2^{1-b}\vare>\vare^b$ since $\vare<\frac12$, and so
	\[x^d\eta^{t}
	=4^{-t}x^{d+d t}
	>\vare^{2t}\vare^{bd(t+1)}
	>\vare^{4bdt}
	=2^{-4bdt\log(1/\vare)}\ge \delta
%	\exp{\left(-10\ell t \log (1/\vare)\right)}
	\]
	since $t\le \frac{5\log(1/\varepsilon)}{\log(\ell(\varepsilon))}$.
	
	Now, let $h:=\abs H$, and let $G$ be such that $\ind_H(G)\le(\delta\abs G)^{h}$. We will show that there exists
$S$ as in the theorem.
	If $\delta\abs G\le1$ then we are done, so we may assume $\delta\abs G>1$, and hence $\abs G>\delta^{-1}>\eta^{-t}$.
	For all $\varepsilon_1,\varepsilon_2\ge\varepsilon$ and every integer $s$ with $0\le s\le t$,
	let $\beta_s(\varepsilon_1,\varepsilon_2)$ be the maximum $\beta>0$ such that for 
	every induced subgraph $F$ of $G$ with $\abs F\ge \eta^s\abs G$, there exists 
	$T\subseteq V(F)$ such that $\abs T\ge \beta\abs F$ and $F[T]$ has edge-density either 
at most $\varepsilon_1$ or at least $1-\varepsilon_2$.
	Since $\delta\le x^{d}\eta^{t}$, it suffices to show that $\beta_0(\varepsilon,\varepsilon)\ge x^{d}\eta^{t}$.
	We claim the following.
	\\
	\\
	(1) {\em For every integer $s$ with $1\le s\le t$,
		and for all $\varepsilon_1,\varepsilon_2\ge\varepsilon$, we have
		\[\beta_{s-1}(\varepsilon_1,\varepsilon_2)
		\ge \eta\cdot\min(\beta_s(p\varepsilon_1,\varepsilon_2),
		\beta_s(\varepsilon_1,p\varepsilon_2)).\]}
\noindent Put $\gamma_1:=\beta_s(p\varepsilon_1,\varepsilon_2)$
		and $\gamma_2:=\beta_s(\varepsilon_1,p\varepsilon_2)$, and let $\gamma=\min(\gamma_1,\gamma_2)$.
		Let $F$ be an induced subgraph of $G$ with $\abs F\ge\eta^{s-1}\abs G$;
we will prove that
                there exists $T\subseteq V(F)$ such that $\abs T\ge \eta\gamma\abs F$ and $F[T]$ has edge-density either
at most $\varepsilon_1$ or at least $1-\varepsilon_2$.
Since $\abs F\ge\eta^{s-1}\abs G$, it follows that 
		$\abs F\ge \eta^{s-1-t}\ge \eta^{-1}$, since $\abs G\ge \eta^{-t}$.
		Since 
		\[\ind_H(F)
		\le \ind_H(G)
		\le (\delta\abs G)^h
		\le (\eta^{-(s-1)}\delta\abs F)^h = (\eta^{t-(s-1)}x^{d}\abs F)^h
		\le (x^{d}\abs F)^h
		\le x^{d}\abs F^h,\]
		there is a blockade $(B_1,\ldots,B_k)$ in $F$, $x$-sparse in one of $F,\overline{F}$,  of length $k\ge \ell(x)$ and width at least 
$$\floor{x^{d}\abs F}=\floor{4\eta\abs F}\ge 2\eta\abs F$$
(since $\abs F\ge \eta^{-1}$).
By the symmetry, we may assume that $(B_1,\ldots,B_k)$ is $x$-sparse in $F$.
		%Let $I$ be a clique or a stable set of $J$ with $q:=\abs I\ge\sqrt{\abs J}\ge \sqrt{\ell(x)}=z^{-1}p$.
		%By the symmetry, we may assume that $I$ is stable in $J$.
		Let $m:=\ceil{\eta\gamma_1\abs F}$.
		
		%Let us renumber $\{B_i:i\in I\}$ as $\{A_1,\ldots,A_q\}$
		%where for all $i,j$ with $1\le i<j\le q$, $A_j$ is $x$-sparse to $A_i$.
		Inductively for $i=k,k-1,\ldots,1$, we choose $C_i\subseteq B_i$ with $|C_i|=m$, as follows.
		Assume $C_k,C_{k-1},\ldots,C_{i+1}$ have been defined,
		and let $D_i$ be their union. 
		Thus $D_i$ is $x$-sparse to $B_i$,
		and so by \ref{tidy} there exists $B_i'\subseteq B_i$ with $\abs{B_i'}\ge\frac12\abs{B_i}$
		such that $B_i'$ is $2x$-sparse to $D_i$; and
		in particular 
$$\abs{B_i'}\ge\frac12\abs{B_i}\ge \eta\abs F \ge \eta^s\abs G.$$
		Thus, by the definition of $\beta_s$,
		there exists $C_i\subseteq B_i'$ with $\abs{C_i}\ge\gamma_1\abs{B_i'}\ge\eta\gamma_1\abs F$
		such that $F[C_i]$ has edge-density either at most $p\varepsilon_1$ or at least $1-\varepsilon_2$.
If its edge-density is at least $1-\varepsilon_2$ then we may set $T=C_i$ and be done; so we may assume that
$F[C_i]$ has edge-density at most $p\varepsilon_1$.
		By averaging, we may assume $\abs{C_i}=\ceil{\eta\gamma_1\abs F}=m$.
		This completes the inductive definition of $C_k,C_{k-1},\ldots,C_1$.
		
		For $1\le i\le k$, $C_i$ is $2x$-sparse to $D_i=C_k\cup C_{k-1}\cup\cdots\cup C_{i+1}$,
		and so there are at most $2xm^2\binom{k}{2}$ edges between $C_1,\ldots,C_k$.
		Therefore, setting $T:=\bigcup_{i=1}^kC_i$, we have 
$$\abs T=km\ge m\ge \eta\gamma_1\abs F\ge \eta\gamma\abs F;$$ 
and
		since $k\ge \ell(x)\ge p/z$ and $2x=(1-z)\varepsilon\le(1-z)\varepsilon_1$, the number of edges of $G[T]$ is at most
		\[kp\varepsilon_1\binom{m}{2}
		+2xm^2\binom{k}{2}
		\le z\varepsilon_1k^2\binom{m}{2}
		+(1-z)\varepsilon_1m^2\binom{k}{2}
		\le \varepsilon_1\binom{km}{2}
		=\varepsilon_1\binom{\abs T}{2}.\]
	So $T$ is a subset of $V(F)$ with the desired property.
	This proves (1).
	
	\bigskip

	By applying (1) for $s=1,2,\ldots,t$, we obtain
	\[\beta_0(\varepsilon,\varepsilon)
	\ge \eta^t\cdot\min_{0\le i\le t}\beta_t(p^{i}\varepsilon,p^{t-i}\varepsilon).\]
	Since $p^t\varepsilon^2\ge 1$ by the choice of $t$, and so $\max\left(p^i\vare, p^{t-i}\vare\right)\ge 1$ if $0\le i\le t$,
	we deduce that $\beta_t(p^i\varepsilon,p^{t-i}\varepsilon)=1$ for all $i$ with $0\le i\le t$;
	and hence $\beta_0(\varepsilon,\varepsilon)\ge \eta^t\ge\delta$, as claimed.
Thus, setting $C=20bd$ satisfies the theorem for all $\varepsilon\in(0,c)$.
	
	Let $a>1$ be such that $2^{-a}=c$;
	we will prove that setting $C=20a^2bd$ works for all $\varepsilon\in(0,\frac12)$.
Thus, let $\varepsilon\in(0,\frac12)$, let 
$$\log(1/\delta)= \frac{20a^2bd\log^2(1/\varepsilon)}{\log(\ell(\varepsilon))},$$
and let $G$ be a graph with $\ind_H(G)\le(\delta\abs G)^{\abs H}$.
Let $\varepsilon':=\varepsilon^a$, and 
let 
$$\log(1/\delta')= \frac{20bd\log^2(1/\varepsilon')}{\log(\ell(\varepsilon'))};$$
then $\varepsilon'\le \varepsilon$, and so $\delta'\ge \delta$
since $\ell$ is non-increasing. Consequently $\ind_H(G)\le(\delta'\abs G)^{\abs H}$. But $\varepsilon'< 2^{-a}=c$, so by what we already proved, 
        there exists
$S\subseteq V(G)$ with $\abs S\ge\delta'\abs G\ge \delta\abs G$ such that one of $G[S],\overline G[S]$ has edge-density at most $\varepsilon'\le \varepsilon$.
	This proves \ref{thm:mainthm2}.~\bbox
	
Finally, let us deduce \ref{mainthm}, which we restate.  Recall that $\mu(G)$
the largest $t$ such that some $t$-vertex induced subgraph of $G$ is a cograph.
\begin{thm}\label{mainthm2}
For every graph $H$ there exists $c>0$ such that $\mu(G)\ge 2^{c\sqrt{\log |G|\log\log |G|}}$ for every $H$-free graph 
$G$ with $|G|\ge 2$.
\end{thm}
\Proof
Let $c$ be as in \ref{mainthmnik}, and let $d=1/(8c)^{1/2}$.
Choose $n_0$ such that for all $n\ge n_0$,
\begin{align*}
\frac14 \log\log n &+ \frac12 \log\log\log n \ge \log (1/d), \text{ and}\\
n^{1/2} & \ge 2^{d\sqrt{\log n\log\log n}}\ge 4.
\end{align*}
Choose $c'\le d/2$ with $c'>0$, such that $n_0^{c'}\le 2$.
We will show that
$\mu(G)\ge 2^{c'\sqrt{\log |G|\log\log |G|}}$ for every $H$-free graph $G$ with $|G|\ge 2$.

Thus, let $G$ be $H$-free. If $|G|\le n_0$, then $|G|^{c'} \le n_0^{c'}\le 2$, and so $\mu(G)\ge |G|^{c'}$. Hence
we may assume that $|G|>n_0$.
Let $\vare=2^{-d\sqrt{\log |G|\log\log |G|}}$, and let
$\delta=2^{-c(\log\frac1\vare)^2/\log\log\frac1\vare}$; then 
\[\log\delta=
-\frac{c(\log\frac1\vare)^2}{\log\log\frac1\vare}
=-\frac{cd^2\log\abs G\log\log\abs G}{\log\log\frac1\vare}.\]
%$$\delta=\exp{\left(\frac{-c\left(\log \frac{1}{\vare}\right)^2}{\log\log \frac{1}{\vare}}\right)}
%= \exp{\left(\frac{-cd^2\log |G|\log\log |G|}{\log\log \frac{1}{\vare}}\right)}.$$
Since 
$$\log\log \frac{1}{\vare}= \frac12 \log\log |G| + \frac12 \log\log\log |G|-\log(1/d) \ge \frac14 \log\log|G|$$
(because $\frac14 \log\log |G| + \frac12 \log\log\log |G| \ge \log (1/d)$), it follows that
$$\log\delta\ge -\frac{cd^2\log |G|\log\log |G|}{\frac14 \log\log |G|}
=-4cd^2\log |G|
=-\frac12\log\abs G,$$
and so $\delta\ge\abs G^{-1/2}$.
By \ref{mainthmnik} and the choice of $c$, there exists $S\subseteq V(G)$ with $|S|\ge |G|^{1/2}$ such that
one of $G[S],\overline G[S]$ has edge-density at most $\vare$. Since $|G|\ge n_0$ it follows that
$|S|>1/\vare$.
By Tur\'an's theorem, $G[S]$ has a clique or stable set of size at least
$$\frac{|S|}{1+\vare|S|}\ge \frac{1}{2\vare}=\frac 12 2^{d\sqrt{\log |G|\log\log |G|}}\ge  
2^{(d/2)\sqrt{\log |G|\log\log |G|}}\ge 2^{c'\sqrt{\log |G|\log\log |G|}}.$$ This proves \ref{mainthm}.~\bbox
%%%%%%%%%%%%%%%%%%%%%%%%%%%%%%%%%%%%%%%%%%%%%%%%%%%%%%%%%%%%%%%%%%%%%%%%%%%%%%%%%%%%%%%%%%%%%%%%%%%%%%%%%%%%%%%%5
\section{Ordered graphs}

An influential paper of Alon, Pach and
Solymosi~\cite{aps} showed that the Erd\H{o}s-Hajnal conjecture admits equivalent
formulations for ordered graphs and for tournaments.  
In this section we observe that our result \ref{mainthm} also extends to ordered graphs and tournaments.
An {\em ordered graph} is a pair $(G,<)$, where $G$ is a graph and $<$ is a linear order of its vertex set.
If $(G,<)$ and $(H,<')$ are ordered graphs, we say $(G,<)$ is {\em $(H,<')$-free} if no induced subgraph of $G$ (made into an ordered graph with the order inherited from $<$ in the natural way)
is isomorphic to $(H,<')$.
Alon, Pach and Solymosi~\cite{aps}
showed that the Erd\H{o}s-Hajnal conjecture \ref{EHconj} is equivalent to the following analogous conjecture for ordered graphs:

\begin{thm}\label{apsconj}
{\bf Conjecture:} For every ordered graph $(H,<')$ there exists $\tau>0$ such that $\kappa(G)\ge |G|^\tau$ for every $(H,<')$-free 
ordered graph $(G,<)$.
\end{thm}
Our theorem \ref{mainthm} translates to:
\begin{thm}\label{orderedmainthm}
For every ordered graph $(H,<')$ there exists $c>0$ such that $\kappa(G)\ge 2^{c\sqrt{\log |G|\log\log |G|}}$ for every $(H,<')$-free 
ordered graph $(G,<)$
with $|G|\ge 2$.
\end{thm}
To prove this, we use 
a theorem of R\"odl and Winkler~\cite{RW}, that says:
\begin{thm}\label{orderedornot} For every ordered graph $(H,<')$, there exists a graph $P$
such that, for
every linear ordering of $V(P)$, the resulting ordered graph is not $(H,<')$-free.
\end{thm}
\noindent{\bf Proof of \ref{orderedmainthm}.\ \ }
Let $(H,<')$ be an ordered graph, and choose $P$ as in \ref{orderedornot}. Choose $c$ satisfying \ref{mainthm} with $H$ replaced by $P$.
Now let $(G,<)$ be an $(H,<')$-free ordered graph. It follows from the property of $P$ that
$G$ is $P$-free, and so by \ref{mainthm}, 
$\kappa(G)\ge 2^{c\sqrt{\log |G|\log\log |G|}}$. This proves \ref{orderedmainthm}.~\bbox

These results also have analogues for tournaments. If $G$ is a tournament, define $\kappa(G)$ to be the size of the largest
transitive subset of $V(G)$.  Our result becomes:
\begin{thm}\label{tourmainthm}
For every tournament $H$ there exists $c>0$ such that $\kappa(G)\ge 2^{c\sqrt{\log |G|\log\log |G|}}$ for every $H$-free 
tournament $G$
with $|G|\ge 2$.
\end{thm}
\Proof
Fix a linear order $<'$ of $V(H)$, and let $H'$ be the graph with vertex set $V(H)$, in which $uv$ is an edge 
if $u$ is earlier than $v$ in the linear order $<'$ and $v$ is adjacent from $u$ in $H$. Thus $(H',<')$ is an ordered graph.
Choose $c$ as in \ref{orderedmainthm} (with $H$ replaced by $H'$). Now let $G$ be an $H$-free tournament.
Derive an ordered graph $(G',<)$ from $G$ similarly. Since $G$ is $H$-free, we deduce that $(G',<)$ is $(H',<')$-free, 
and the result follows
from \ref{orderedmainthm}, since every clique or stable set of $G'$ is a transitive set of $G$. This proves \ref{tourmainthm}.~\bbox

\end{document}